\documentclass[12pt,reqno]{amsart}
\usepackage{amsmath, amsfonts, amssymb, amsthm, hyperref}
\usepackage{bm}
\allowdisplaybreaks[4]
\def\l{\left}
\def\r{\right}
\def\bg{\bigg}
\def\({\bg(}
\def\){\bg)}
\def\t{\text}
\def\f{\frac}

\def\ls{\leqslant}

\def\da{\delta}
\def\la{\lambda}

\def\Proof{\noindent{\it Proof}}

\def\Z{\mathbb Z}
\def\C{\mathbb C}

\def\Q{\mathbb Q}

\def\1{{\bf 1}}

\def\<{\langle}
\def\>{\rangle}
\theoremstyle{plain}
\newtheorem{theorem}{Theorem}[section]
\newtheorem{lemma}{Lemma}
\newtheorem{corollary}{Corollary}

\newtheorem{proposition}{Proposition}

\theoremstyle{definition}

\theoremstyle{remark}
\newtheorem{remark}{Remark}

\newcommand{\sign}[1]{\mathrm{sign}(#1)}
\makeatletter
\@namedef{subjclassname@2020}{%
  \textup{2020} Mathematics Subject Classification}
\makeatother
 \vspace{4mm}

\begin{document}
\hbox{Preprint,{\tt arXiv:2206.02589}}
\medskip

\title[Proof of a conjecture involving derangements and roots of unity]
{Proof of a conjecture involving
\\ derangements and roots of unity}
\author{Han Wang}
\address {(Han Wang) Department of Mathematics, Nanjing
University, Nanjing 210093, People's Republic of China}
\email{hWang@smail.nju.edu.cn}

\author{Zhi-Wei Sun}
\address {(Zhi-Wei Sun, corresponding author) Department of Mathematics, Nanjing
University, Nanjing 210093, People's Republic of China}
\email{zwsun@nju.edu.cn}

\keywords{Derangement, Hermitian matrix, determinant, roots of unity, eigenvalue.
\newline \indent Supported by the Natural Science Foundation of China (grant no. 11971222).}
\begin{abstract}
Let $n>1$ be an odd integer, and let $\zeta$ be a primitive $n$th root of unity in the complex field.
Via the Eigenvector-eigenvalue Identity, we show that
$$\sum_{\tau\in D(n-1)}\mathrm{sign}(\tau)\prod_{j=1}^{n-1}\frac{1+\zeta^{j-\tau(j)}}{1-\zeta^{j-\tau(j)}}
=(-1)^{\frac{n-1}{2}}\frac{((n-2)!!)^2}{n},
$$
where $D(n-1)$ is the set of all derangements of $1,\ldots,n-1$.
This confirms a previous conjecture of Z.-W. Sun. Moreover, for each $\delta=0,1$ we determine the value of $\det[x+m_{jk}]_{1\ls j,k\ls n}$ completely, where
$$m_{jk}=\begin{cases}(1+\zeta^{j-k})/(1-\zeta^{j-k})&\text{if}\ j\not=k,\\\delta&\text{if}\ j=k.
\end{cases}$$
\end{abstract}
\maketitle

\section{Introduction}
\setcounter{lemma}{0}
\setcounter{theorem}{0}
\setcounter{equation}{0}
\setcounter{conjecture}{0}
\setcounter{remark}{0}
\setcounter{corollary}{0}

For $n\in\Z^+=\{1,2,3,\ldots\}$, let $S_n$ be the symmetric group of all permutations of $\{1,\ldots,n\}$. A permutation $\tau\in S_n$ is called a {\it derangement} of $1,\ldots,n$ if $\tau(j)\not=j$
for all $j=1,\ldots,n$. For convenience, we use $D(n)$ to denote the set of all derangements of $1,\ldots,n$. The derangement number $D_n=|D(n)|$ plays important roles in enumerative combinatorics.
It is well known that
$$D_n=n!\sum_{k=0}^n\f{(-1)^k}{k!}$$
(cf. (10.2) of \cite[p.\,90]{vW}).

Let $n>1$ be an odd integer. Z.-W. Sun \cite[Theorem 1.2]{S19} proved that
$$\det\l[\tan\pi\f{j-k}n\r]_{1\ls j,k\ls n-1}=n^{n-2}.$$
As
$$\tan\pi x=\f{2\sin\pi x}{2\cos\pi x}=i\f{1-e^{2\pi ix}}{1+e^{2\pi ix}},$$
we see that
\begin{align*}\det\l[\tan\pi\f{j-k}n\r]_{1\ls j,k\ls n-1}&=i^{n-1}\det\l[\f{1-\zeta^{j-k}}{1+\zeta^{j-k}}\r]_{1\ls j,k\ls n-1}
\\&=(-1)^{(n-1)/2}\sum_{\tau\in D(n-1)}\sign{\tau}\prod_{j=1}^{n-1}\f{1-\zeta^{j-\tau(j)}}{1+\zeta^{j-\tau(j)}},
\end{align*}
where $\zeta=e^{2\pi i/n}$.

Z.-W. Sun \cite{S21} and \cite[Conj. 11.24]{S-book} conjectured that if $n>1$ is odd and $\zeta$
is a primitive $n$th root of unity in the complex firld $\C$ then
\begin{equation}\label{1}
\sum_{\tau\in D(n-1)}\sign\tau\prod_{j=1}^{n-1}\f{1+\zeta^{j-\tau(j)}}{1-\zeta^{j-\tau(j)}}
=(-1)^{\f{n-1}{2}}\f{((n-2)!!)^2}{n},
\end{equation}
Our first goal is to prove an extension of this conjecture.

\begin{theorem}\label{Th1.1} Let $n>1$ be an odd integer, and let $\zeta\in\C$ be a primitive $n$th root of unity. For $j,k=1,\ldots,n$ define
$$a_{jk}=\begin{cases}(1+\zeta^{j-k})/(1-\zeta^{j-k})&\t{if}\ j\not=k,
\\0&\t{if}\ j=k.\end{cases}$$
Then we have
\begin{equation}\label{x-form}
\det[x+a_{jk}]_{1\ls j,k\ls n-1}=(-1)^{\f{n-1}{2}}\f{((n-2)!!)^2}{n}.
\end{equation}
\end{theorem}

Applying Theorem \ref{Th1.1} with $x=1$, we immediately obtain the following result.
\begin{corollary}
Let $n>1$ be odd. Then, for any primitive $n$th root $\zeta\in\C$ of unity, we have
\[
\det[\tilde a_{jk}]_{1\ls j,k\ls n-1}=(-1)^{\f{n-1}{2}}\f{((n-2)!!)^2}{n2^{n-1}},
\]
where
\[
\tilde a_{jk}=\begin{cases}
1/(1-\zeta^{j-k})&\t{if}\ j\neq k,\\
1/2&\t{if}\ j=k.\end{cases}
\]
\end{corollary}

For any odd integer $n>1$, Sun \cite{S21} also conjectured that
if $\zeta\in\C$ is a primitive $n$th root of unity then
\begin{equation}\label{Guo}
\sum_{\tau\in D(n-1)}\sign\tau\prod_{j=1}^{n-1}\f{1}{1-\zeta^{j-\tau(j)}}=\f{(-1)^{\f{n-1}{2}}}n\l(\f{n-1}2!\r)^2.
\end{equation}
Quite recently, X. Guo et al. \cite{Guo} proved \eqref{Guo} via using the following result
which dates back to Jacobi in 1834 (cf. P.B. Denton, S.J. Parke, T. Tao and X. Zhang \cite[Theorem 1]{Tao}).

\begin{theorem} [Eigenvector-eigenvalue Identity]\label{Th-Tao} Let $A$ be an $n\times n$ matrix over $\C$ which is Hermitian (i.e., the transpose $A^T$ of $A$ coincides with the conjugate of $A$), and let $\la_1,\ldots,\la_n$
be its $n$ real eigenvalues. Let $v_n=(v_{n,1},\ldots,v_{n,n})^T$ be an eigenvector associated with the eigenvalue $\la_n$ of the matrix $A$ such that its norm $\|v_n\|=\sqrt{\sum_{j=1}^n|v_{n,j}|^2}$
equals $1$.
Let $j\in\{1,\ldots,n\}$ and let $A_j$ be the $(n-1)\times(n-1)$ Hermitian matrix formed by deleting the $j$th row and the $j$th column from A. Let $\la_{j,1},\ldots,\la_{j,n-1}$
 be all the real eigenvalues of $A_j$. Then we have
 $$|v_{n,j}|^2\prod_{k=1}^{n-1}(\la_n-\la_k)=\prod_{k=1}^{n-1}(\la_n-\la_{j,k}).$$
\end{theorem}

Motivated by Theorem \ref{Th1.1}, we also establish the following result.

\begin{theorem}\label{Th1.3}
Let $n>1$ be odd. Then,  for any primitive $n$th root $\zeta\in\C$ of unity, we have
\begin{equation}\label{xb}
\det[x+b_{jk}]_{1\ls j,k\ls n-1}=(-1)^{\f{n+1}2}(nx+1)\f{((n-1)!!)^2}{n(n-1)},
\end{equation}
where
\[
b_{jk}=
\begin{cases}
(1+\zeta^{j-k})/(1-\zeta^{j-k})&\t{if}\ j\not=k,\\
1&\t{if}\ j=k.\end{cases}
\]
\end{theorem}

 We are going to prove Theorems \ref{Th1.1} and \ref{Th1.3} in Sections 2 and 3 respectively.

\section{Proof of Theorem \ref{1}}
\setcounter{lemma}{0}
\setcounter{theorem}{0}
\setcounter{equation}{0}
\setcounter{conjecture}{0}
\setcounter{remark}{0}
\setcounter{corollary}{0}

We need the following easy lemma.

\begin{lemma}\label{Lem2.1} Let $n\in\Z^+$ and $s\in\{0,\ldots,n-1\}$. For any primitive $n$th root $\zeta$
of unity in a field $F$, we have the identity
\begin{equation}\label{x-iden}\sum_{0<r<n}\f{\zeta^{-rs}}{1-x\zeta^r}=\f{\sum_{j=0}^{n-1}x^j-nx^s}{x^n-1}.
\end{equation}
\end{lemma}
\Proof. Clearly,
\[
\sum_{r=0}^{n-1}\f{\zeta^{-rs}}{1-x\zeta^r}=\sum_{r=0}^{n-1}\f{\zeta^{-rs}}{1-x^n}
\sum_{k=0}^{n-1}(x\zeta^r)^k=\sum_{k=0}^{n-1}\f{x^k}{1-x^n}\sum_{r=0}^{n-1}\zeta^{r(k-s)}=\f{nx^s}{1-x^n}.
\]
Thus
\[
\sum_{r=1}^{n-1}\f{\zeta^{-rs}}{1-x\zeta^r}=\f{nx^s}{1-x^n}-\f{1}{1-x}=\f{\sum_{j=0}^{n-1}x^j-nx^s}{x^n-1}
\]
as desired. \qed

\begin{remark} Lemma \ref{Lem2.1} in the case $F=\C$ is essentially equivalent to \cite[Theorem 3.1]{Gessel}.
\end{remark}

\begin{corollary} \label{+-1} Let $n\in\Z^+$ and $s\in\{0,\ldots,n-1\}$. Let $\zeta$ be any primitive $n$th root of unity in the field $\C$.

{\rm (i)} If $n$ is odd, then
\begin{equation}\label{+1}\sum_{0<r<n}\f{\zeta^{-rs}}{1+\zeta^r}=\f{(-1)^sn-1}2.
\end{equation}

{\rm (ii)} We have
\begin{equation}\label{-1}\sum_{0<r<n}\f{\zeta^{-rs}}{1-\zeta^r}=\f{n-1}2-s.
\end{equation}
\end{corollary}
\Proof. (i) When $n$ is odd, putting $x=-1$ in \eqref{x-iden} we immediately get \eqref{+1}.

(ii) Letting $x\to1$ in \eqref{x-iden} we obtain \eqref{-1} since
\begin{align*}\lim_{x\to1}\f{\sum_{j=0}^{n-1}x^j-nx^s}{x^n-1}
&=\lim_{x\to1}\f{(\sum_{j=0}^{n-1}x^j-nx^s)'}{(x^n-1)'}=\lim_{x\to1}\f{\sum_{0<j<n}jx^{j-1}-nsx^{s-1}}
{nx^{n-1}}
\\&=\f{\sum_{j=0}^{n-1}j-ns}{n}=\f{1}{n}\sum_{j=0}^{n-1}j-s=\f{n-1}{2}-s.
\end{align*}
by L'Hospital's rule.

Combining the above, we have completed the proof of Corollary \ref{+-1}. \qed

\begin{remark} It seems that the identity \eqref{-1} should be known long time ago.
We note that it essentially appeared as \cite[(3.5)]{Gessel} though $(n-1)/2$
in \cite[(3.5)]{Gessel} should be corrected as $(n+1)/2$.
\end{remark}

Now we give an auxiliary proposition.

\begin{proposition}\label{Prop} Let $n\in\Z^+$, $k\in\{1,\ldots,n\}$ and $s\in\{0,\ldots,n-1\}$.
For any primitive $n$th root $\zeta$ of unity in a field $F$, we have
\begin{equation}\label{gen}\sum_{j=1\atop j\not=k}^n\f{1+x\zeta^{j-k}}{1-x\zeta^{j-k}}\zeta^{s(k-j)}
=1+2\f{\sum_{j=0}^{n-1}x^j-nx^s}{x^n-1}-n\da_{s,0}.
\end{equation}
Consequently, if $\zeta$ is a primitive $n$th root of unity in $\C$, then
\begin{equation}\label{2}
\sum_{j=1\atop j\neq k}^{n}\f{1+\zeta^{j-k}}{1-\zeta^{j-k}}\zeta^{s(k-j)}=
\begin{cases}
n-2s&\t{if}\ 0<s<n,\\
0&\t{if}\ s=0.\end{cases}
\end{equation}
\end{proposition}
\Proof. In view of Lemma \ref{Lem2.1}, we have
\begin{align*}
\sum_{j=1\atop j\neq k}^{n}\f{1+x\zeta^{j-k}}{1-x\zeta^{j-k}}\zeta^{s(k-j)}
&=\sum_{r=1}^{n-1}\f{1+x\zeta^r}{1-x\zeta^r}\zeta^{-sr}
=2\sum_{r=1}^{n-1}\f{\zeta^{-rs}}{1-x\zeta^r}-\sum_{r=1}^{n-1}\zeta^{-rs}
\\&=2\f{\sum_{j=0}^{n-1}x^j-nx^s}{x^n-1}+1-\sum_{r=0}^{n-1}\zeta^{-rs}
=2\f{\sum_{j=0}^{n-1}x^j-nx^s}{x^n-1}+1-n\da_{s,0}.
\end{align*}
This proves \eqref{gen}.

When $F=\C$, letting $x\to1$ in \eqref{gen} or using the identity \eqref{-1}, we get \eqref{2}.
\qed

We also need another lemma.

\begin{lemma}[Sun \cite{S19}] \label{MM'}
For any matrix $M=[m_{jk}]_{0\ls j,k\ls n}$ over $\C$, we have
\[
\det[x+m_{jk}]_{0\ls j,k\ls n}=\det(M)+x\det(M'),
\]
where $M'=|m_{jk}'|_{1\ls j,k\ls n}$ with $m_{jk}'=m_{jk}-m_{j0}-m_{0k}+m_{00}$.
\end{lemma}

\medskip
\noindent{\it Proof of Theorem \ref{Th1.1}}.
Obviously $A=[a_{kj}]_{1\ls k,j\ls n}$ is a Hermitian matrix.
For each $k=1,\ldots,n$, by Proposition \ref{Prop} we have
$$\sum_{j=1}^na_{kj}\zeta^{-js}=\sum_{j=1\atop j\not=k}^n\f{1+\zeta^{j-k}}{1-\zeta^{j-k}}\zeta^{-js}
=\begin{cases}(n-2s)\zeta^{-ks}&\t{if}\ s\in\{1,\ldots,n-1\},
\\0&\t{if}\ s=n.\end{cases}$$
Thus $\la_s=n-2s\ (s=1,\ldots,n-1)$ and $\la_n=0$ are all the eigenvalues of $A$.
Moreover, for each $s=1,\ldots,n$, the column vector $$v^{(s)}=\f1{\sqrt n}(\zeta^{-s},\zeta^{-2s},\ldots,\zeta^{-ns})^T$$ is an eigenvector
of norm $1$ associated with the eigenvalue $\la_s$.

Let $A_n$ be the Hermitian matrix $[a_{kj}]_{1\ls k,j\ls n-1}$, and let
$\la_{n,1},\ldots,\la_{n,n-1}$ be all the eigenvalues of $A_n$.
Note that $v^{(n)}=(1,\ldots,1)^T/\sqrt n$. Applying Theorem \ref{Th-Tao} with $j=n$, we obtain that
$$(-1)^{n-1}\det(A_n)=\prod_{k=1}^{n-1}(0-\la_{n,k})=\bg|\f1{\sqrt n}\bg|^2\,\prod_{k=1}^{n-1}(0-\la_k)=\f{(-1)^{n-1}}n\prod_{k=1}^{n-1}(n-2k)
$$
and hence
$$\det(A_n)=\f1n\prod_{k=1}^{(n-1)/2}(n-2k)(n-2(n-k))=\f{(-1)^{(n-1)/2}}n\prod_{k=1}^{(n-1)/2}(n-2k)^2
=\f{(-1)^{(n-1)/2}}n((n-2)!!)^2.$$
On the other hand,
\begin{equation*}\det(A_n)=\det(A_n^T)=\sum_{\tau\in D(n-1)}\sign{\tau}\prod_{j=1}^{n-1}\f{1+\zeta^{j-\tau(j)}}{1-\zeta^{j-\tau(j)}}.
\end{equation*}
Combining the last two equalities, we immediately get \eqref{x-form} for $x=0$.

By Lemma \ref{MM'}, we have
$$\det[x+a_{jk}]_{1\ls j,k\ls n-1}=\det(A_n)+x\det(A_n'),$$
where $A_n'=[a_{jk}']_{2\ls j,k\ls n-1}$ with
$$a_{jk}'=a_{jk}-a_{j1}-a_{1k}+a_{11}=a_{jk}-a_{j1}-a_{1k}.$$
It is easy to see that $a_{kj}'=-a_{jk}'$ for all $j,k=2,\ldots,n-1$.
So we have
$$\det(A_n')=\det(-A_n')=(-1)^{n-2}\det(A_n')=-\det(A_n')$$
and hence
$$\det[x+a_{jk}]_{1\ls j,k\ls n-1}=\det(A_n)+x\det(A_n')=\det(A_n)=\f{(-1)^{(n-1)/2}}n((n-2)!!)^2.$$
This ends our proof. \qed

\section{Proof of Theorem \ref{Th1.3}}
\setcounter{lemma}{0}
\setcounter{theorem}{0}
\setcounter{equation}{0}
\setcounter{conjecture}{0}
\setcounter{remark}{0}
\setcounter{corollary}{0}

\begin{lemma} \label {Lem-c} Let $n\in\{2,3,4,\ldots\}$, and let $\zeta$ be a primitive $n$th root of unity.
For $j,k=1,\ldots,n$ define
$$
c_{jk}=\begin{cases}
1/(1-\zeta^{j-k})&\t{if}\ j\ne k,\\
0&\t{if}\ j=k.\end{cases}$$

{\rm (i)} The $n$ eigenvalues of $[c_{jk}+\da_{jk}]_{1\ls j,k\ls n}$ are $s-\f{n-1}2\ (s=1,\ldots,n)$.

{\rm (ii)} If $n$ is odd, then
\begin{equation}\det[c_{jk}+\da_{jk}]_{1\ls j,k\ls n-1}=(-1)^{\f{n+1}2}\f{(n+1)((n-1)!!)^2}{n(n-1)2^{n-1}}.
\end{equation}
\end{lemma}
\Proof. (i) For $j,k=1,\ldots,n$ let
$$t_{jk}=\begin{cases}1+i\cot\pi\f{j-k}n&\t{if}\ j\not=k,
\\0&\t{if}\ j=k.\end{cases}.$$
By F. Calogero and A. M. Perelomov \cite[Theorem 1]{CP}, the $n$ numbers $2s-n-1\ (s=1,\ldots,n)$
are all the eigenvalues of the matrix $[t_{jk}]_{1\ls j,k\ls n}$.
Thus
\begin{equation}\label{mjk}\det[xI_n-t_{jk}]_{1\ls k\ls n}=\prod_{s=1}^n(x-(2s-n-1)),
\end{equation}
where $I_n$ be the identity matrix of order $n$.
For $j,k=1,\ldots,n$ with $j\not=k$, clearly
$$t_{jk}=1-\f{2\cos\pi\f{j-k}n}{2i\sin\pi\f{j-k}n}=
1-\f{e^{2\pi i\f{j-k}n}+1}{e^{2\pi i\f{j-k}n}-1}=\f2{1-e^{2\pi i\f{j-k}n}}.$$
Note that $\zeta=e^{2\pi ia/n}$ for some $1\ls a\ls n$ with $\gcd(a,n)=1$.
 Applying the Galois automorphism $\sigma_a$ in the Galois group $\mathrm{Gal}(\Q(e^{2\pi i/n})/\Q)$
 with $\sigma_a(e^{2\pi i/n})=e^{2\pi ia/n}$, we obtain from \eqref{mjk} the polynomial identity
\begin{equation}\label{CP-Iden}
\det[xI_n-2c_{jk}]_{1\ls k\ls n}=\prod_{s=1}^n(x-(2s-n-1)).
\end{equation}
Thus
\[
\det[xI_n-c_{jk}|]_{1\ls j,k\ls n}=\prod_{s=1}^n\l(x-s+\f{n+1}2\r),
\]
and hence
\begin{align*}
\det[xI_n-c_{jk}-\da_{jk}]_{1\ls j,k\ls n}=&\det[(x-1)I_n-c_{jk}]_{1\ls j,k\ls n}
\\=&\prod_{s=1}^n\l(x-1-s+\f{n+1}2\r)=\prod_{s=1}^n\l(x-\l(s-\f{n-1}2\r)\r).
\end{align*}
So the numbers $s-\f{n-1}2\ (s=1,\ldots,n)$ are all the eigenvalues of $[c_{jk}+\da_{jk}]_{1\ls j,k\ls n}$.

(ii) Now assume that $n$ is odd.
Let $$\{\la_1,\ldots,\la_n\}=\l\{\f{3-n}2, \f{5-n}2,\ldots,\f{n+1}2\r\}$$ with $\la_n=0$. Then the column vector $$v^{(n)}=\f1{\sqrt n}(\zeta^{-\f{n-1}2},\zeta^{-2\f{n-1}2},\ldots,\zeta^{-n\f{n-1}2})^T$$ is an eigenvector
of norm $1$ associated with the eigenvalue $\la_n$.

Let $C_n$ be the Hermitian matrix $[c_{kj}+\da_{jk}]_{1\ls k,j\ls n-1}$, and let
$\la_{n,1},\ldots,\la_{n,n-1}$ be all the eigenvalues of $C_n$.
Note that $v^{(n)}=(\zeta^{-\f{n-1}2},\ldots,\zeta^{-\f{n(n-1)}2})^T/\sqrt n$. Applying Theorem \ref{Th-Tao} with $j=n$, we obtain that
$$(-1)^{n-1}\det(C_n)=\prod_{k=1}^{n-1}(0-\la_{n,k})=\bg|\f{\zeta^{-\f{n(n-1)}2}}{\sqrt n}\bg|^2\,\prod_{k=1}^{n-1}(0-\la_k)=\f{(-1)^{n-1}}n\prod_{k=1 \atop k\ne \f{n-1}{2}}^{n}\l(k-\f{n-1}2\r)
$$
and hence
\[\begin{aligned}
\det(C_n)=&\f{(n-1)(n+1)}{2^{n-1}n}\prod_{k=1}^{(n-3)/2}(n-1-2k)(n-1-2(n-1-k))\\
=&(-1)^{(n+1)/2}\f{(n-1)(n+1)}{2^{n-1}n}\prod_{k=1}^{(n-3)/2}(n-1-2k)^2
=(-1)^{(n+1)/2}\f{(n+1)((n-1)!!)^2}{2^{n-1}n(n-1)}.\end{aligned}
\]
This concludes the proof. \qed

\medskip
\noindent{\it Proof of Theorem \ref{Th1.3}}.
Let $B$ be the $n\times n$ matrix $[b_{kj}]_{1\ls k,j\ls n}$.
With the aid of \eqref{2},
\begin{equation}\label{3}
1+\sum_{j=1\atop j\neq k}^{n}\f{1+\zeta^{j-k}}{1-\zeta^{j-k}}\zeta^{s(k-j)}=
\begin{cases}
n+1-2s&\t{if}\ 0<s<n,\\
1&\t{if}\ s=0.\end{cases}
\end{equation}
Thus, for each $k=1,\ldots,n$, we have
$$\sum_{j=1}^nb_{kj}\zeta^{-js}=\zeta^{-ks}+\sum_{j=1\atop j\not=k}^n\f{1+\zeta^{j-k}}{1-\zeta^{j-k}}\zeta^{-js}
=\begin{cases}(n+1-2s)\zeta^{-ks}&\t{if}\ s\in\{1,\ldots,n-1\},
\\1&\t{if}\ s=n.\end{cases}$$

Recall that $n$ is odd. Let
$$\{\mu_1,\ldots,\mu_n\}=\{n-1, n-3,\ldots,2,1,0,-2,\ldots,-n+3\}$$ with $\mu_n=0$. Then the column vector $$u^{(n)}=\f1{\sqrt n}(\zeta^{-\f{n+1}2},\zeta^{-2\f{n+1}2},\ldots,\zeta^{-n\f{n+1}2})^T$$ is an eigenvector
of norm $1$ associated with the eigenvalue $\mu_n$.

Let $B_n$ be the Hermitian matrix $[b_{kj}]_{1\ls k,j\ls n-1}$, and let
$\mu_{n,1},\ldots,\mu_{n,n-1}$ be all the eigenvalues of $B_n$.
Note that $u^{(n)}=(\zeta^{-\f{n+1}2},\ldots,\zeta^{-\f{n(n+1)}2})^T/\sqrt n$. Applying Theorem \ref{Th-Tao} with $j=n$, we obtain that
$$(-1)^{n-1}\det(B_n)=\prod_{k=1}^{n-1}(0-\mu_{n,k})=\bg|\f{\zeta^{-\f{n(n+1)}2}}{\sqrt n}\bg|^2\,\prod_{k=1}^{n-1}(0-\mu_k)=\f{(-1)^{n-1}}n\prod_{k=1\atop k\ne \f{n+1}{2}}^{n-1}(n+1-2k)
$$
and hence
\[\begin{aligned}
\det(B_n)=&\f{n-1}n\prod_{k=2}^{(n-1)/2}(n+1-2k)(n+1-2(n+1-k))\\
=&(-1)^{(n+1)/2}\f{n-1}n\prod_{k=1}^{(n-1)/2}(n+1-2k)^2
=(-1)^{(n+1)/2}\f{((n-1)!!)^2}{n(n-1)}.\end{aligned}
\]
This proves \eqref{xb} for $x=0$.

By Lemma \ref{MM'} we have
$$\det[x+b_{jk}]_{1\ls j,k\ls n-1}=\det(B_n)+x\det(B_n')$$
for certain $(n-2)\times(n-2)$ matrix $B_n'$ over $\C$ not depending on $x$.
As $1+b_{jk}=2(c_{jk}+\da_{jk})$ (with $c_{jk}$ given by Lemma \ref{Lem-c}) for all $j,k=1,\ldots,n-1$,
we have
\begin{align*}\det(B_n)+\det(B_n')
=&\det[1+b_{jk}]_{1\ls j,k\ls n}=2^{n-1}\det[c_{jk}+\da_{jk}]_{1\ls j,k\ls n-1}
\\=&(n+1)(-1)^{(n+1)/2}\f{((n-1)!!)^2}{n(n-1)}=(n+1)\det(B_n)
\end{align*}
with the aid of Lemma \ref{Lem-c}. Therefore
\begin{align*}\det[x+b_{jk}]_{1\ls j,k\ls n-1}=&\det(B_n)+x(n\det(B_n))=(1+nx)\det(B_n)
\\=&(-1)^{(n+1)/2}(1+nx)\f{((n-1)!!)^2}{n(n-1)}
\end{align*}
as desired. This ends our proof of Theorem \ref{Th1.3}. \qed

\end{document}